\documentclass[a4paper]{amsart}

\usepackage{amsmath, amssymb,latexsym, amsthm, amsfonts, amscd}
\usepackage{verbatim}
\usepackage{epsfig}
\usepackage{graphics}
\usepackage{citesort}

\title[Classification of Tangle Solutions for Integrases]{Classification of Tangle
Solutions for Integrases, A Protein Family that Changes {D}{N}{A}
Topology}

\author[D.~Buck]{Dorothy~Buck}
\address{Department of Mathematics, Imperial College London}
\email{d.buck@imperial.ac.uk}

\author[C.~Verjovsky~Marcotte]{Cynthia~Verjovsky~Marcotte}
\address{Department of Mathematics, St. Edward's University}
\email{cynthm@admin.stedwards.edu}

\newtheorem{thm}{Theorem}[section]
\newtheorem{THM*}{Theorem}
\newtheorem{cor}[thm]{Corollary}
\newtheorem{lem}[thm]{Lemma}

\theoremstyle{Definition}

\theoremstyle{remark}

\numberwithin{equation}{section}

\newcommand{\nibf}{\noindent \textbf}
\newcommand{\niit}{\noindent \textit}

\newcommand{\oc}{O_c}
\newcommand{\of}{O_f^k}

\newcommand{\otf}{O_f^2}
\newcommand{\onotf}{O_f^{k \neq 2}}

\newcommand{\Qonly}{\mathbb{Q!}}
\newcommand{\Q}{\mathbb{Q}}

\newcommand{\Z}{\mathbb{Z}}

\newcommand{\wt}{\widetilde}

\newcommand{\beq}{\begin{equation}}
\newcommand{\eeq}{\end{equation}}
\newcommand{\bea}{\begin{eqnarray*}}
\newcommand{\eea}{\end{eqnarray*}}

\newcommand{\dbc}{{\rm{dbc}}}
\newcommand{\core}{{\rm{core}}}

\begin{document}
\maketitle
\section*{\today}




%

\section*{Abstract}

Integrase proteins acting on circular double-stranded DNA often
change its topology by transforming unknotted circles into torus
knots and links. Two systems of tangle equations---corresponding
to the two initial DNA sequences---arise when modelling this
transformation: direct and inverted.

With no \textit{a priori} assumptions on the constituent tangles,
we utilize Dehn surgery arguments to completely classify the
tangle solutions for each of the two systems. A key step is to
combine work of our previous paper \cite{us} with recent results
of Kronheimer, Mrowka, Ozsv\'{a}th and Szab\'{o}~\cite{KMOS} and
Ernst~\cite{E1} to show a certain prime tangle must in fact be a
generalized Montesinos tangle.

These tangle solutions are divided into three classes, common to
both systems, plus a fourth class for the inverted system that
contains the sole generalized Montesinos tangle. We discuss the
possible biological implications of our classification, and of
this novel solution.

\section{Introduction}

DNA is often circular, and so can be knotted or linked
\cite{LiuWang}. To aid in essential cellular tasks, many proteins
have thus evolved to manipulate the geometry and topology of DNA.
Important examples of these tasks include: \emph{replication}
(creating 2 copies from one double strand of DNA),
\emph{inversion} (inverting a subsequence of DNA), and
\emph{integration} and \emph{deletion} (inserting and deleting DNA
sequences into/from other DNA segments). For instance, certain
proteins invert crossings of two double strands of DNA to unknot
or unlink DNA \cite{BergWang}.

It is important for biologists to understand how these proteins
operate. So a model is often generated based on data obtained by
experimentally probing a particular system. For example, if a
particular protein can knot or link circular DNA, then by running
DNA through a gel, one can determine the crossing number of the
DNA knots or links. One can then use electron microscopy to
identify the exact knot or link type~\cite{KrasStas}.

By modelling different regions of the (possibly knotted) DNA
molecule as tangles, we can describe the protein's action as a
change in one of the constituent tangles. Mathematicians can help
biologists by finding all tangle combinations that may explain the
protein's action. Some of these mathematically possible solutions
can then be eliminated through biological considerations.

Based on the biological work of Wasserman and
Cozzarelli~\cite{WDC,WC}, Ernst and Sumners developed the
tangle model to describe and make predictions---later
experimentally verified---about how the protein Tn3 interacts with
DNA \cite{ES1}. The tangle model has since been used to determine
specific features of a wide variety of protein-DNA interactions
(see for
example,~\cite{us,Cri,ES2,E2,Dar,Dar3,SakaVaz,DLV,me,SECS,WhiMilCozz,Vaz}).

\subsection{The integrase family of recombinases}
We will focus on one particular family of proteins that affects
DNA topology: the integrase family of site-specific recombinases,
also known as the tyrosine recombinases. For an overview of this
family, see the review of Grainge and Jayaram~\cite{GJ}. The
integrase family includes \emph{Flp}, $\lambda$ \emph{Int}, and
\emph{Cre}. Integrases are involved in a wide variety of
biological activities, including integrating (\textit{i.e.},
inserting) a virus' DNA into a host cell's DNA. Integrase proteins
have a common biochemical mechanism and share a subsequence of
amino acids. Given varying numbers of DNA axis self-crossings
(\textit{supercoils}) these proteins can transform unknotted
circular DNA into a variety of torus knots or links~\cite{GJ}.

Members of the other family of site-specific recombinases, the
serine recombinases, act very differently. These proteins, such as
Tn3, require a fixed number of supercoils. The protein then cuts
and rearranges the DNA (occasionally multiple times) before
releasing it. For several serine recombinases, it has been
possible to completely solve the tangle model equations
(see~\cite{ES1,ES2,Vaz}).

In contrast, for integrase proteins, the unbound DNA's
configuration varies from molecule to molecule, depending on the
degree of supercoiling. Further, these proteins act only once
before releasing the DNA. So, to model the varying configurations,
and its effects on the resulting DNA products after a single round
of recombination, we must use a larger number of tangles than
needed for a serine recombinase. Rather than four fixed tangles, a
solution to the (direct or inverted) system of equations consists
of three fixed tangles, $P$, $R$ and $\oc$, and an indexed family
of tangles, $\of$. In general, there are an infinite number of
solutions to such a system. The increased complexity of the tangle
model for the integrase family has thus far prevented a full
solution to the tangle model. Previous work found solutions by
making several assumptions which are thought to be biologically
reasonable and which lead to considerable mathematical
simplification~\cite{Cri,SECS,me}.

We previously considered the tangle model for a generic
member of the integrase family~\cite{us}
without any \textit{a priori} assumptions on the constituent
tangles. We gave and exemplified three
(topologically motivated) classes of solutions for $\of, \oc$ and
$P$, all of which are valid for both systems. We then developed Dehn
Surgery arguments to eliminate all (in the direct system) or all
but one (in the inverted system) other possibilities.

In this paper, we now complete the classification of all possible
solutions to the systems of equations arising from both the direct
and inverted cases. In particular, we utilize recent work of
Kronheimer, Mrowka, Ozsv\'{a}th and Szab\'{o}~\cite{KMOS}, and
Ernst~\cite{E1,E2} to determine the fourth solution in the
inverted case. In addition, we determine all possible solutions
for $R$ in both the direct and inverted settings. Combining the
results derived below with our previous work~\cite{us}, we obtain
a complete classification of all solutions to the equations
arising from either of the direct or inverted system:

\begin{THM*}
There are three classes of solutions to the equations that model
integrase-mediated DNA recombination
\begin{tabbing}
\label{invdown}
Before:\ \= \kill
\>$N(\of + \oc + P) = b(1,1) =$ unknot for $k \in \{0,1,2,3\}$\\
\>$N(\of + \oc + R) = b(2k,1) =$ {$T_{2k,2}$ torus link} for
direct sites\\
\> \hspace{2in}\textbf{or}\\
\>$N(\of + \oc + R) = b(2k+1,1) =${$T_{2k+1,2}$ torus knot} for
inverted sites\\

\end{tabbing}

\begin{enumerate}
\item[1.] $P= (\infty)$, $\oc$, $\of$ and $R$
are integral tangles.
%
\item[2.] $P$ and $\oc$ are integral, $\of = (\infty)$ for at most one value of $k$,
integral for at most two values of $k$ and otherwise is vertical or
the sum of a vertical plus an integral tangle.
$R = (\infty)$, integral or vertical plus
integral. However if $\of=(\infty)$ for some $k$, then $R$ cannot be
$(\infty)$, nor can it be integral unless $k=0$.

\item[3.] $P$ and $R$ are strictly rational,  $\oc$ is integral, and $\of$
is integral for at most value of $k$ and otherwise strictly
rational.
\end{enumerate}

\vspace{0.05in}

\noindent For the direct system, there are no other
classes of solutions.

\noindent For the inverted system, there is precisely one additional
class of solutions:
\begin{enumerate}
\item[4.] \
$P = \pm(p)$, $R = \pm(1+p)$, $\oc$ is integral, $\otf$ is the
prime Montesinos tangle $\mp(\frac{1}{2},\frac{2}{3},p-1)$ and
$\onotf = \pm(-p+\frac{1}{2k})$ or $\mp(p+\frac{1}{2k+2})$, where
$p\in\{0,1\}$.
\end{enumerate}
%
%
\end{THM*}

\vspace{-0.05in}

Theorem~1 is a consequence of Theorem~1 in~\cite{us} and
Theorem~\ref{finalinvsoln}, as well as
Theorems~\ref{R1},~\ref{Rvert} and~\ref{ThirdSol}. Note that $\oc$
is always integral. Also note that the only thing that keeps the
solutions in Class 4 from being solutions in Class 2 is tangle
$\otf$.

Note that solutions are unique up to compensating negative and
positive vertical sums of vertical tangles or horizontal sums of
horizontal tangles. See Section 2 for a discussion of these
equivalence classes.

The paper is organized as follows: In Section 2 we review some
basic facts about tangles, four-plats and their corresponding
double branch covers. In Section 3 we provide the biological
motivation and background for our work: the action of an integrase
on DNA with either inverted or direct repeats. We also present the
generalization of Ernst and Sumners' mathematical model in terms
of tangles and four-plats~\cite{ES1}. In Section 4 we focus on the
inverted system. Here we harness recent results of Kronheimer,
Mrowka, Ozsv\'{a}th and Szab\'{o}~\cite{KMOS} and
Ernst~\cite{E1,E2} to limit the possibilites for $\otf$, and then
determine the complete fourth solution.

In Section~\ref{Solns} we present examples of all classes of
solutions to these tangle equations, including the novel inverted
solution. We conclude with some remarks on the biological
relevance of our model and solutions, as well as possible
directions for future (biological and mathematical) work.

\section{Tangles, Four-Plats and their Double Branch Covers}

We begin by recalling a few elementary facts about tangles. (For a
more comprehensive introduction to tangles, see~\cite{us} and
\cite{GK}). A \textit{tangle} $T$ is a pair $(B^3,t)$, where $B^3$
is a 3-ball with a given boundary parametrization with four
distinguished boundary points labelled NW, NE, SW, SE, and $t$
consists of a pair of properly embedded unoriented arcs with
endpoints NW, NE, SW and SE. We say two tangles $A$ and $B$ are
{\textit{equivalent}} if there exists an isotopy taking $A$ to
$B$, which remains the identity on $\partial A$.

Tangles can be divided into three mutually exclusive families:
locally knotted, rational and prime.

A tangle is \textit{locally knotted} if there exists a sphere in $B^3$
meeting $t$ transversely in 2 points such that the 2-ball bounded
by the sphere intersects $t$ in a knotted spanning arc. Locally knotted
tangles have reducible double branch covers (branching over $t$).

We note that locally knotted tangles do not occur in our context
since the DNA molecule starts off unknotted, $P$, $\oc$ and $\of$
cannot be locally knotted. Suppose $R$ were locally knotted, with
$S\subset R$ be a sphere bounding a ball $B$ containing the
locally knotted arc. Then, since every product is a prime knot,
$N(\of+\oc+R)-B$ must be an unknotted arc, for all $k$. But this
implies that $N(\of+\oc+R)=N(O^j_f+\oc+R)$ for $j\neq k$, a
contradiction. So none of the constituent tangles are locally
knotted.

\textit{Rational} tangles are the second family; they are so
called because their equivalence classes are in one-to-one
correspondence with the extended rational numbers
($\Q\cup\{\infty\}$) via a continued fraction expansion, as first
constructed by Conway~\cite{Con}. (See~\cite{GK,KL1} for nice
classifications.) A tangle whose corresponding rational number is
$\frac{p}{q}$ will be denoted by $\left(\frac{p}{q}\right)$.
Rational tangles are formed by an alternating series of horizontal
and vertical half-twists of two (initially untwisted) parallel
arcs (and hence are freely isotopic to them). Any continued
fraction decomposition of $\frac{p}{q} = a_n + 1/(a_{n-1}+ ... \
(1/a_1))$
yields a finite list of integers $[a_1,\ldots,a_n]$ which tell us
how to twist the strands around each other to get a diagram of the
tangle. The $(0)$ tangle corresponds to two untwisted horizontal
arcs (one joining NE to NW and one joining SE to SW), whereas the
$(\infty)$ tangle corresponds to two untwisted vertical arcs.
The double branch cover of a rational tangle is a solid torus.

All locally unknotted, non-rational tangles are \textit{prime}.
The double branch cover branched over $t$ of a prime tangle is
irreducible and has incompressible boundary~\cite{Lic}. Bleiler
demonstrated that the minimal prime tangle has a minimal
projection with five crossings~\cite{Ble}.

There are several operations one can perform on tangles. We
concentrate on three. The first operation forms a knot or
2-component link from a given tangle $A$: the \emph{numerator
closure}, $N(A)$. This adds an unknotted arc joining the northern
endpoints, and another unknotted arc joining the southern
endpoints, or equivalently, the boundary of $A$ and $(0)$ are
identified so that $\textrm{E}_A$ is identified with
$\textrm{E}_{(0)}$ for $E\in\{NE,NW,SE,SW\}$.

The second operation, \emph{tangle sum}, takes a pair of tangles
$A$, $B$, and, under certain restrictions, yields a third tangle,
$A + B$, by identifying the eastern hemispheric boundary disk of
$A$ with the western one of $B$ in such a way that $\textrm{NE}_A$
is identified with $\textrm{NW}_B$ and $\textrm{SE}_A$ is
identified with $\textrm{SW}_B$. Note that the (0) tangle is the
identity under this operation: $A+(0)=A$. Beware that under tangle
addition, we cannot distinguish between $A+(p)$ added to $B+(-p)$
and $A$ added to $B$. Thus although tangle summands are written in
their simplest form, they are unique only up to an arbitrary
number of compensating positive and negative horizontal twists.

The last operation is the \emph{vertical sum} $\star$, which takes
two tangles $A$ and $B$ and, under certain restrictions, yields a
third tangle, $A\star B$, by identifying the southern hemispheric
boundary disk of $A$ with the northern hemispheric boundary disk
of $B$ in such a way that $\text{SE}_A$ is identified with
$\text{NE}_B$ and $\text{SW}_A$ is identified with $\text{NW}_B$.
Note that under the numerator closure operation, given two tangles
$A$ and $B$, although $(A\star (1/n) )+ B\not\simeq A+ (B\star
(1/n))$ in fact$N\left((A\star (1/n) )+ B\right)\simeq N\left(A+
(B\star (1/n))\right)$. Also, $N\left((A\star (1/n) )+ (B\star
(-1/n))\right)\simeq N\left(A+ B\right)$. As with horizontal
tangle sums, we cannot distinguish between these two cases. Thus
although tangle solutions are written in the simplest form for a
given equivalence class, they are unique only up to an arbitrary
number of compensating positive and negative vertical twists.

A particular class of prime tangles is obtained by tangle sum of
rational tangles and will be of interest to us: \emph{Montesinos}
tangles~\cite{Mon}. We will use the notation
$\left(\frac{a_1}{b_1},\frac{a_2}{b_2},\ldots\frac{a_n}{b_n}\right)$
to denote the Montesinos tangle obtained by the tangle sum of the
rational (possibly integral) tangles
$\left(\frac{a_1}{b_1}\right)$,
$\left(\frac{a_2}{b_2}\right),\ldots
\left(\frac{a_n}{b_n}\right)$. A tangle that is ambient isotopic
to a sum of rational tangles called \emph{a generalized Montesinos
tangle}. Since the sum of a rational tangle and an integral tangle
yields a rational tangle~\cite{Qua}, a Montesinos (non-rational)
tangle must have at least two non-integral summands.

The numerator closure of a rational tangle yields a
\textit{four-plat}, a knot or 2-component link that admits a
projection consisting of a braid on 4 strings, with one strand
free of crossings~\cite{BS}. Given two rational tangles
$\left(\frac{p}{q}\right)$ and $\left(\frac{p'}{q'}\right)$ then
$N\left(\frac{p}{q}\right) = N\left(\frac{p'}{q'}\right)$ iff $p =
p'$ and $q^{\pm1} \equiv q \mod p$. Thus given a 4-plat, we can
write it as the numerator closure of a rational tangle that is
unique up to the relationship above. (See~\cite{KL2} for a
classification of rational knots.) Schubert showed that all
four-plats are prime knots~\cite{Sch}. A four-plat obtained by the
numerator closure of $\left(\frac{p}{q}\right)$ is written as
$b(p,q)$. For example, we can write the unknot as $b(1,1)$, and
the trefoil as $b(3,1)$.


\subsection{Double Branch Covers}

If $T$ is a tangle, then $\wt{T}$ will mean the double cover of
$B^3$, branched over $t$. In general, we will write $\dbc(K)$ to
denote the three-manifold that is the double cover of $S^3$
branched over the set $K$. We now turn our attention to the
(compact, connected and orientable) three-manifolds that arise as
double branch covers of tangles or four-plats.

If $P$ is a rational tangle, then $\wt{P}$ is a solid  torus,
which we will denote by $V_P$.
Schubert showed that $\dbc(b(p,q))$ is the lens space $L(p,q)$.
Two four-plats $b(p,q)$ and $b(p',q')$ are equivalent if and only
if their corresponding double branch covers, the lens spaces
$L(p,q)$ and $L(p',q')$, are homeomorphic~\cite{Sch}, so
$b(p,q)=b(p',q')$ if and only if $p=p'$ and $q^{\pm 1}\equiv q'
\mod p$. (See Rolfsen~\cite{Rolf}, for a good introduction to lens
spaces.)

A summing disk of a tangle (either the western or eastern summing
one) lifts to an annulus on the boundary of the double branch
cover. This annulus can be \emph{meridional} (its core bounds a
meridional disk of the solid torus), \emph{longitudinal} (its core
intersects a meridian once) or neither (intersects a meridian more
than once). When two tangles are summed, this corresponds to
gluing two solid tori along the annuli that are lifts of the
summing disks. Hence, a natural subdivision of rational tangles
arises: We say a rational tangle is \textit{integral}, and write
it as $(n)$, if it consists of a series of $n$ horizontal
half-twists, where $n \in \mathbb{Z}$. We denote this class as
$\mathbb{Z}$. Integral tangles have summing disks that lift to
longitudinal annuli. Similarly, a tangle is the infinity tangle,
denoted by $(\infty)$, if it consists of two vertical parallel
strands with no twists. The infinity tangle has a summing disk
that lifts to a meridional annulus. We say a tangle is
\textit{strictly rational}, and denote this class as $\Qonly$ if
it is neither integral nor the infinity tangle, and so $\mathbb{Q}
= \Qonly \cup \mathbb{Z} \cup \{ (\infty)\}$. Strictly rational
tangles have summing disks that lift to annuli that are neither
meridional nor longitudinal. Within the class of strictly rational
tangles we distinguish \textit{vertical} tangles, written as
$\left(\frac{1}{n}\right)$, which consists of a series of $n$
vertical half-twists ($|n|>1$).

\section{Biological Motivation and Model}

We can now describe in full detail a generalization~\cite{SECS} of
the original tangle model of Ernst and Sumners~\cite{ES1}.

We illustrate our model with a member the integrase family of
recombinases, the protein \textit{Flp} (pronounced `flip'). Flp
has served as the paradigm for site-specific recombination
\cite{YangJay}, and there are a number of proteins (including Cre,
and $\lambda$ Int acting on LR sites) whose products are,
topologically speaking, identical to those of Flp.

Roughly speaking, Flp recognizes two copies of a specific DNA
sequence, binds at these sites, cuts the DNA at the sites, moves
the strands in space, reseals the break, and releases the DNA.
When acting on circular DNA, Flp can change the underlying knot
type of the DNA, for example turning the unknot into the trefoil
knot. The distribution of knot/link products should reflect the
supercoiling density \cite{PollNash}. We call a DNA molecule that
has not been acted on by Flp a \textit{substrate}, and a molecule
that has been acted on a \textit{product}. In these terms, the
substrate is always an unknot and the products are torus knots or
links.

We model each of the substrates and products as the numerator
closure of the sum of three tangles. Each tangle arc represents a
segment of double-stranded DNA. In the tangle model pioneered by
Ernst and Sumners~\cite{ES1}, the cutting and joining of DNA is
assumed to be completely localized: two of the tangles are
unchanged by the action of the protein. In the substrate, the
first tangle, $P$ ($P$arental), represents two short identical
sites that Flp recognizes and to which it chemically binds and
then cuts, rearranges and re-seals. This action can be thought of
as removing $P$ and replacing it with a new tangle, $R$
($R$ecombinant), in the product. The second tangle, $\oc$,
represents the part of the DNA that is physically constrained, but
unchanged, by the protein ($O$ stands for Outside and $c$ for
constrained). The last tangle, $\of$, represents the part of the
DNA that is free (hence the subscript $f$) from protein binding
constraints. $\of$ can vary depending on the amount of DNA
supercoiling present at the time Flp acts. The superscript $k$
indexes these different possibilities.

In terms of tangles, this amounts to saying that our substrate and
products can be modelled as:

\begin{tabbing}
this is as long\= \kill
\>$N(\of + \oc + P) =$ \textrm{substrate (before recombination)}  \\
\>$N(\of + \oc + R) =$ \textrm{product (after recombination)}{ {
\rule{0cm}{0.5cm}}}
\end{tabbing}
where $k \in \{0,1,2,3\}$. $\of$ varies as $k$ varies, so we
obtain different products, as described below. We use $O^k$ to
mean the part unchanged by Flp, that is, $O^k = \of + \oc$. When
there is no mathematical distinction between $\of$ and $\oc$, we
use $O_1$ and $O_2$ to represent them interchangeably,
\emph{i.e.}, $\{O_1,O_2\}=\{\of,\oc\}$.

Recall that, in contrast to integrases, proteins in the serine family of
recombinases, such as Tn3, require a fixed number of supercoils
before they begin cutting and rejoining DNA. Once this requirement
is met, they rearrange the DNA, occasionally multiple times,
before releasing it. The corresponding tangle equations: from
substrate $N(O_f+O_c+P)=K_0$ to products $N(O_f+O_c+nR)=K_n$,
were first solved (\textit{i.e.}, all constituent tangles have
been characterized, given the 4-plats) by Ernst and
Sumners~\cite{ES1}. Note that the free part, $O_f$, does not vary.
This single, fixed $O_f$ is what has thus far made the serine
recombinase tangle equations more tractable than the integrases,
whose equations involve a family of tangles $\of$, indexed by $k$.

\subsection{Two Systems:  Inverted and Direct}

Flp identifies two short identical sequences, called
\textit{repeats}, on a molecule of DNA. These sites are
non-palindromic sequences, and can thus be given an orientation,
and hence on circular DNA, the strings can be in head to head
(\textit{inverted repeats}) or head to tail (\textit{direct
repeats}) orientation. Action on inverted repeats on a circular
molecule of DNA yields a knot, and action on direct repeats yields
a two-component link. When Flp acts on DNA it yields a variety of
torus knots (inverted repeats) and links (direct repeats) that
depend on $\of$.

When Flp acts on a DNA molecule with inverted sites, experiments
have shown that the resulting DNA can be an unknot (with a
different DNA sequence), or a knot with up to 11 crossings
\cite{GJ}. Crisona \emph{et al.} have obtained images (using
electron microscopy) of the simplest products, and has shown that
they are the torus knots $b(\pm 1,1)$ (the unknot), primarily
positive $b(3,1)$ and exclusively positive $b(5,1)$~\cite{Cri}.
This experimental evidence indicates that Flp begins with an
unknotted DNA substrate with inverted repeats, $b(\pm 1,1)$ and
converts it via tangle surgery into a torus knot $b(\pm
(2k+1),1)$, where $k \in \{0,1,2,3\}$. (The chirality of the
products for $k = 1$ and $k=2$ has not been determined for all
members of the integrase family, so we remain in the general
situation. Corollary \ref{chiral} considers the specific setting
for Flp (and $\lambda$ Int acting on LR sites), whose products'
handedness are known.)

We thus model the action of Flp on DNA with inverted repeats as:
\smallskip

\begin{tabbing}
\label{invdown}
Before:\ \= \kill
\textrm{Before:} \>$N(\of + \oc + P) = b(\pm 1,1) =$ unknot, for \ $k \in \{0,1,2,3\}$\\
\textrm{After:} \>$N(O^0_f + \oc  + R) = b(\pm 1,1) =$ unknot { { \rule{0cm}{0.5cm}}}\\
\>$N(O^k_f + \oc + R) = b(\pm (2k+1),1) =$ torus knot for \ $k \in
\{0,1,2,3\}${{ \rule{0cm}{0.5cm}}}
\end{tabbing}

\smallskip

When Flp acts on a DNA molecule with direct sites, experiments
have shown that the resulting DNA can be an unlink, or a
2-component link with up to 10 crossings~\cite{GJ}.
Electrophoretic gels have determined that the simplest products
are $b(0,1)$, $b(\pm 2,1)$ and $b(\pm 4,1)$~\cite{GJ}. This
experimental evidence indicates that Flp begins with an unknotted
DNA substrate with direct repeats, $b(\pm 1,1)$ and converts it
via tangle surgery into a torus link $b(\pm 2k,1)$, where $k \in
\{0,1,2,3\}$. We thus model the action of Flp on DNA with direct
repeats as:

\smallskip

\begin{tabbing}
\label{invdown}
Before:\ \= \kill
Before: \>$N(\of + \oc + P) = b(\pm 1,1) =$ unknot, \textrm{for}\ $k \in \{0,1,2,3\}$ \\
After:\>$N(O^0_f + \oc + R) = b(0,1) =$ unlink { { \rule{0cm}{0.5cm}}}\\
\>$N(O^k_f + \oc + R) = b(\pm 2k,1) =$ torus link for \ $k \in
\{0,1,2,3\}${{ \rule{0cm}{0.5cm}}}
\end{tabbing}
\smallskip

\subsection{Strategy}

Given the set of tangle equations above, whose products (4-plats)
are known, the goal is to determine the constituent tangles. The
interplay of tangles and four-plats with their corresponding
double branch covers is the key to many of our results in tangle
calculus. For instance, if $C$ and $D$ are tangles, and $D$ is a
rational tangle, then the $\dbc(C+D)$ is obtained by gluing
$\wt{C}$ and $\wt{D}=V_D$ along annuli that are the lifts of their
corresponding gluing disks. If $D$ is integral, then the gluing
annulus is boundary reducible, and $\dbc(C+D)\simeq \wt{C}$.

The sum and subsequent numerator closure of two tangles $C$ and
$D$ induces a gluing of the boundaries of their respective double
branch covers $\wt{C}$ and $\wt{D}$. If $N(C+D)$ yields a
four-plat $b(p,q)$, then $\wt{C} \cup_h \wt{D}$ must be the lens
space $L(p,q)$, where $h$ is the map that takes $\mu_{\partial
\wt{C}}$ to $p\lambda_{\partial \wt{D}} + q\mu_{\partial \wt{D}}$.
In particular, when $C$ and $D$ are both rational, $\wt{C}=V_C$
and $\wt{D}=V_D$ are solid tori, and they form a Heegaard
splitting $V_C\cup_h V_D$ of $L(p,q)$.

Replacing tangle $P$ in $N(O+P)$ by tangle $R$ to obtain $N(O+R)$
is called \emph{tangle surgery}. If $P$ and $R$ are rational
tangles, then tangle surgery corresponds to replacing $V_P$ with
$V_R$ in the double branch cover, and thus corresponds to
different Dehn fillings of $\wt{O^k}$. In the case of
$N(O^k+P)=b(1,0)$, the unknot, then since $\dbc(b(1,0))$is $S^3$,
and the tangle surgery corresponds to Dehn surgery on a knot
complement ($O^k$) in $S^3$. If $O^k$ is not rational, then the knot
is non-trivial.

Previous work (Ernst and Sumners~\cite{ES1}, and the authors
\cite{us} via different techniques for $R$) proved that $P$ and $R$
are rational for both direct and inverted repeats. So the tangle
surgery of replacing $P$ with $R$ corresponds to Dehn surgery in
the double branch covers.

Thus the strategy is to use restrictions on the type of Dehn
surgeries of $S^3=\dbc(b(1,1))$ that yield lens spaces. This in
turn restricts the possible tangle solutions.

\section{The fourth case for inverted repeats}

In~\cite{us}, we asked whether there were any solutions in the
single remaining open inverted case: $P$ rational, $O_1$ integral,
$O_2$ prime (and hence $O^k = \of + \oc$ prime). We now can give a
positive answer for $k=2$ and eliminate this case for all other
$k$.

In this section, we consider only the inverted system, as we
previously completely classified the direct system~\cite{us}.

\subsection{Double branch cover of
\mathversion{bold}$\otf$\mathversion{normal} is a trefoil knot complement}

We begin by restricting the possibilities for $O^k$, in part by
harnessing a powerful recent result of Kronheimer, Mrowka,
Ozsv\'{a}th and Szab\'{o}:

\begin{thm}[\cite{KMOS}, Corollary 8.4]
\label{KMOS} If $K$ is a knot in $S^3$, such that for some
$r\in\mathbb{Z}$, $M_K(r)=L(p,q)$ is a lens space where $|p | <9$,
then $K$ must either be the unknot or the trefoil knot.
\end{thm}

\begin{thm}
In the inverted repeats system, $O^k$ is rational for $k\in\{ 0,
1, 3\}$. Further, if $O^2$ is not rational, then $\wt{O^2}$ must
be the complement of a trefoil knot.~\label{trefoil}
\end{thm}

\nibf{Proof:} Recall that since both $P$ and $R$ are rational
tangles, then the tangle equations $\ N(O^k + P) = b(1,1)$ and $\
N(O^k + R) = b(2k+1,1)$ correspond to a Dehn surgery along
$\core(V_P)$ in the (possibly trivial)
knot complement $\wt{O^k} := S^3 \setminus
V_P$ that yields $L(2k+1,1)$.

The possibilities of $\wt{O^k}$ being satellite (for all $k$) or
torus (for $k \neq 2$) have been eliminated in Theorem 7.1
of~\cite{us}. Gordon~\cite{Gor} and Moser~\cite{Mos} have
classified all surgeries on a generic torus knot $T_{a,b}$
complement that yield lens spaces: $L(p,qb^2)$ iff $p = qab \pm
1$. A straightforward calculation shows that only $L(5,4)$ can be
obtained from a torus knot---the trefoil knot $T_{3,2}$.

It now remains only to rule out the possibility of
$\wt{O^k}$ being hyperbolic.
Assume $\wt{O^k}$ is a hyperbolic knot ($K$) complement. Then by
the Cyclic Surgery Theorem~\cite{CGLS1}, the surgery slope must be
integral, and in fact must be $2k+1$, since
$H_1(M_K(s/t))=\Z_{|s|}$. In our setting, the integral surgery
slope $2k+1$ is strictly less than $9$, since $k\in\{0,1,2,3\}$.
Hence, we can apply Corollary 8.4
of~\cite{KMOS}, and obtain that
$\wt{O^k}$ must be a solid torus or a trefoil knot
complement, neither of which is hyperbolic. Therefore the
hyperbolic case is impossible.

We have thus shown that $\wt{O^{k \neq 2}}$ must be the complement
of the unknot (a solid torus), and therefore $O^k$ is a rational
tangle. Further, $\wt{O^2}$ is either a solid torus or trefoil
knot complement, and therefore $O^2$ is either a rational tangle
or a tangle whose double branch cover is a trefoil knot
complement. \hfill $\Box$

\subsection{\mathversion{bold}$\oc$\mathversion{normal} is integral
and \mathversion{bold}$\otf$\mathversion{normal} is prime}

We can now use the rationality of $P$ and $O^k := \oc + \of$ for
$k \neq 2$ to place restrictions on the summand tangles.

\begin{thm}
In the inverted system, assume $O^2 = O_1 + O_2$ has a trefoil
knot complement double branch cover. Then $\oc$ must be integral
and $\otf$ must be prime with $\wt{\otf}$ is a trefoil knot
complement.~\label{trefoil2}
\end{thm}

\nibf{Proof:}
%
We first show that our only option in this setting is $O_1$ is
integral and $O_2$ is prime.

From Section 6 (as summarized in Table 1) of our previous paper
\cite{us}, the only cases not eliminated are ($i$) $P$ rational,
$O_1$ prime and $O_2$ integral, and ($ii$) $P$ integral, $O_1 =
(\infty)$ and $O_2$ prime. We now show possibility ($ii$) cannot
occur.

Theorem 6.8($v$) of~\cite{us}
eliminates the case when $\oc = O_1 = (\infty)$ and $\of
= O_2$ is prime. Alternately, the case with $\oc$ prime and $\of =
(\infty)$ can occur for at most 1 value of $k$, since for $k \neq
j$, $N(\of + \oc + R) \neq N(O_f^j + \oc + R)$. By Theorem
6.8($viii$) in~\cite{us}, $\of$ must be integral for all
other values of $k$. But then for these 3 (or 4) values of $k$,
$O^k := \of + \oc = $ prime + integral, is prime~\cite{Qua}, a
contradiction to Theorem~\ref{trefoil}.

Thus the only possibility that can occur is $O_1$ is integral and
$O_2$ is prime. Suppose that $\oc = O_2$ is prime and $\otf = O_1$
is integral. Then since $\oc+\of$ is rational for $k\neq 2$, then
$\of = (\infty)$ for $k = 0$, $1$ and $3$ by Cam Van~\cite{Qua}.
But for $k \neq j$, $O^k \neq O^j$ as these are different
recombination products. Thus $\oc$ must be integral, and $\otf$
prime.

Since $\oc$ is integral, $\dbc(\oc+\otf)\simeq\wt{\otf}\simeq
\wt{O^2}$, the trefoil knot complement from the preceding theorem.
\hfill $\Box$

\subsection{\mathversion{bold}$\otf$\mathversion{normal} is a Montesinos tangle,
and \mathversion{bold}$P$\mathversion{normal} and
\mathversion{bold}$R$\mathversion{normal} are integral}

Montesinos links were first considered by Bonahon~\cite{Bon} and
Montesinos~\cite{Mon}, by using work of Tollefson~\cite{Toll} who
determined that every involution of a Seifert fiber space with
non-empty boundary must respect a Seifert fibration. Therefore if
the double branch cover of a link in $S^3$ admits a Seifert
fibration which is invariant under the covering involution the
link is either a torus link or what is now called a Montesinos
link. Ernst~\cite{E1}, with a later clarification by Darcy
\cite{Dar}, used related techniques to determine that any tangle
whose double branch covers is Seifert fiber space must be a
rational or a generalized Montesinos tangle:

\begin{thm}[Ernst~\cite{E1}, and Darcy \cite{Dar}]\label{montesinos}
If M is a SFS with orbit surface a disk and $n \geq 0$ exceptional
fibers and if M is the 2-fold branch cover of a tangle (B,t), then
B is a generalized Montesinos tangle.
\end{thm}

\begin{cor} In the inverted system, the
case where $P$ is a rational tangle, $\oc$ is an integral tangle,
and $\otf$ is a prime tangle, then $\otf$ is a generalized
Montesinos tangle with two non-integral rational tangles.
~\label{montesinos-sol}
\end{cor}

\nibf{Proof.} From Theorem~\ref{trefoil2} we have that $\wt{\otf}$
is a trefoil knot complement. By Theorem~\ref{montesinos}, $\otf$
then must be a generalized Montesinos tangle with two non-integral
summands.~\hfill$\Box$

\smallskip

The following result was first used without proof by Darcy in
\cite{Dar}. We give a short self-contained proof below, since we
will also need this for the fourth solution.

\begin{cor}
If $N(\otf+0)$ is a four-plat, then $\otf$ is of the form
$\left(\frac{a}{b},\frac{c}{d}\right)\star\left(\frac{1}{m}\right)$
with $m \ \epsilon \ \mathbb{Z}$. \label{newcor}
\end{cor}

\nibf{Proof.} From Corollary 4.5, we know that $\otf$ is a
generalized Montesinos tangle with two non-integral rational
tangles. Note that
$\left(\frac{a}{b},\frac{c}{d}\right)+(m)=\left(\frac{a}{b},\frac{c+dm}{d}\right)$,
a Montesinos tangle with two non-integral summands.

As a consequence of the above and \cite{Dar2}, such a generalized
Montesinos tangle is of the form
\[
\left(\left(\left.\left(\frac{a}{b},\frac{c}{d}\right)\star\left(\frac{1}{m_n}\right)\right)+(m_{n-1})\right)\star
\ldots+(m_1)\right)
\]
if $n$ is even, or
\[
\left(\left(\left.\left(\frac{a}{b},\frac{c}{d}\right)\star\left(\frac{1}{m_n}\right)\right)+(m_{n-1})\right)\star
\ldots\star\left(\frac{1}{m_1}\right)\right)
\]
if $n$ is odd.

We will show that if $n> 1$ then $N(\otf+0)\simeq N(\otf)$ cannot
be a four-plat. Note first that if $n$ is odd, then
$N\left(\left(\left.\left.\left(\frac{a}{b},\frac{c}{d}\right)\star\left(\frac{1}{m_n}\right)\right)+(m_{n-1})\right)
\star
\ldots+\left(m_{2}\right)\right)\star\left(\frac{1}{m_1}\right)\right)\simeq
 N\left(\left(\left(\left(\frac{a}{b},\frac{c}{d}\right)\star\left(\frac{1}{m_n}\right)\right)+(m_{n-1})\right)\star
\ldots+\left(m_{2}\right)\right)$. Hence we need only examine the
case when $n>1$ is even.

Recall that $N\left((A\star(1/n))+B\right)\simeq
N\left(A+(B\star(1/n))\right)$. Hence by a simple inductive
argument we can show that
$N\left(\left(\left(\left(\frac{a}{b},\frac{c}{d}\right)\star\left(\frac{1}{m_n}\right)\right)+(m_{n-1})\right)\star
\ldots+\left(m_{1}\right)\right)\simeq
N\left(\left(\frac{a}{b},\frac{c}{d}\right)
+\left(\left((m_{n-1})\ldots+\left(\left((m_3)
+\left((m_1)\star\left(\frac{1}{m_2}\right)\right)\right)\star\left(\frac{1}{m_4}\right)
\ldots\right)\star\left(\frac{1}{m_n}\right)\right)\right.\right)$.
The tangle $\left(\left((m_{n-1})\ldots+\left(\left((m_3)
+\left((m_1)\star\left(\frac{1}{m_2}\right)\right)\right)\star\left(\frac{1}{m_4}\right)
\ldots\right)\star\left(\frac{1}{m_n}\right)\right)\right.$ is a
rational tangle, isotopic to the tangle given by the vector $[m_1,
m_2, \ldots, m_n,0]$. If $n>1$, this is not an integral tangle,
and so
 $N\left(\left(\left.\left(\frac{a}{b},\frac{c}{d}\right)\star\left(\frac{1}{m_n}\right)\right)+(m_{n-1})\right)\star
\ldots+(m_1)\right)$ is isotopic to the numerator closure of the
Montesinos tangle with \emph{three} non-integral rational tangles
whose double-branch cover is a Seifert fiber space with three
exceptional fibers, and thus cannot be a four-plat. Thus $n=1$,
and so $\otf$ is of the form
$\left(\frac{a}{b},\frac{c}{d}\right)\star\left(\frac{1}{m}\right)$.
Note if $m=1$, then $\otf$ is properly (not generalized)
Montesinos. \hfill$\Box$

\begin{thm}[Ernst~\cite{E2}]
%
%
In the inverted system, if $O^2$ a Montesinos tangle with two
non-integral rational summands, then $P$ and $R$ are integral
tangles.~\label{PRint}
\end{thm}

\nibf{Proof.}
Recall that $N(O^2+P)=b(1,1)$, and $N(O^2+R)=b(5,1)$, whose double branch
covers are $S^3$ and $L(5,1)$ respectively. Now $O^2$ is a
Montesinos tangle with two non-integral summands. If $P$ were not
integral, then the double branch cover of $N(O^2+P)$ would be a
Seifert fibered space with three exceptional fibers, a
contradiction, as it is a lens space. Similarly with $R$.
\hfill$\Box$

\subsection{The Final Inverted Solution}

From the previous section, this fourth case reduces to $P$, $R$
and $\oc$ are integral and $\otf$ is Montesinos with two
non-integral tangles. We now determine all solutions to the tangle
equations in this case, by applying an algorithm of Ernst whose
completeness is guaranteed~\cite[Theorem 3.1]{E2}.

\begin{thm}
In the fourth solution of the inverted case, the tangle solutions
are $P = \pm(p)$, $\oc$ integral, $R = \pm(1+p)$, $\otf =
\mp(\frac{1}{2},\frac{2}{3},p-1)$ and $\onotf =
\pm(-p+\frac{1}{2k})$ or $\mp(p+\frac{1}{2k+2})$, where
$p\in\{0,1\}$.

\label{finalinvsoln}
\end{thm}

\nibf{Proof:} First note that if $P=(p)$, where $p\neq 0$, we can
assume that $P=(0)$ by moving the $p$ horizontal twists into
$O^k$. So given a solution $\hat{P}=(0)$, $\hat{R}=(r)$, and
$\hat{O}^k$, then the solution corresponding to $P=(p)$ would be
$R=(r+p)$ and $O^k=\hat{O}^k+(-p)$.

\nibf{Case \mathversion{bold}1: $k=2$\mathversion{normal}.} From
Theorems~\ref{trefoil} and~\ref{PRint}, $O^2$ is a Montesinos
tangle of two non-integral summands, and $P$ and $R$ are integral
tangles. In this section we will use the same notation as in the
algorithm in~\cite{E2}: $O^2=(u/v,x/y)$, $P=(m_0)$ and $R=(m)$. By
moving the $m_0$ horizontal twists of $P$ into $O^2$, we can set
$P = (0)$, $R = (m)$, and $O^2 = (u/v,x/y,m_0)$. We can now apply
the algorithm to determine $m_0$, $m$, $u$, $v$, $x$ and $y$ .

For $k=2$, we begin with $b(\alpha_1,
\beta_1) = b(1,0)$ and obtain $b(\alpha_2,\beta_2) = b(+5,1)$.
Ernst's algorithm yields $m_0 = m = -1$, $u = 1$, $v = x = 2$ and
$y = 3$. So the solution is either $R=(-1)$, and $O^2 =
(\frac{u}{v},\frac{x}{y},m_0) = (\frac{1}{2},\frac{2}{3},-1) =
(\frac{1}{2}, \frac{-1}{3})$ or their mirror images: $R=(1)$, and
$O^2 = (-\frac{1}{2},-\frac{2}{3},1)$.

\nibf{Case \mathversion{bold}2: $k\neq 2$\mathversion{normal}.} We
will show that if $P=(0)$, then $O^k=\mp\left(\frac{1}{2k}\right)$ or
$\pm\left(\frac{1}{2k+2}\right)$.

We have that $N(O^k+P)=N(O^k+(0))=N(O^k)=
N\left(\frac{1}{r}\right)$ and
$N(O^k+R)=N(O^k\mp(1))=N\left(\pm\frac{2k+1}{1}\right)$.

By the correspondence between numerator closure of rational tangles and
 4-plats,
$N\left(\frac{p}{q}\right)=N\left(\frac{p'}{q'}\right)$ if and
only if $p=p'$ and $q^{\pm1}\equiv q'\mod
p$~\cite{Con},~\cite{Sch}. This tells us that $O^k$ is a vertical
tangle $\left(\frac{1}{q}\right)$. The tangle
$(O^k+R)=\left(\frac{1}{q}\pm 1\right)= \left(\frac{1\pm
q}{q}\right)$, where the sign is the sign of $R$. Hence $1\pm q=
\pm 2k+1$ and so $q=2k$ or $-2k-2$ when $R=(1)$, and $q=-2k$ or
$2k+2$ when $R=(-1)$.

\smallskip

\noindent We thus obtain the general solution by converting as
described above.
 \hfill $\Box$

Note that as $\oc$ is integral, then $\oc$ is either
zero or the integral part of $\of$.

We reiterate that, although solutions are given in their simplest
form, they are unique only up to compensating positive and
negative vertical twists. So for example, as noted by Darcy, $P =
(0)$, $R = (1)$ and $O^2 = (\frac{1}{2},\frac{2}{3}) \star
(\frac{-1}{2})$ is also a solution for $k=2$ \cite{Dar4}.

For both Flp and $\lambda$ Int (acting on LR sites), Crisona
\textit{et al.} have numerous electron microscope images of two of
the products~\cite{Cri}. Most of trefoils ($k=1$) and all of the
pentafoils ($k=2$) examined are positive. This chirality further
constrains the tangle possibilities for Flp and $\lambda$ Int (on
LR sites).

\begin{cor}\label{chiral}
For both Flp and $\lambda$ Int (acting on LR sites), if $P=(p)$,
then $p\in\{0,-1\}$ and $R = (p-1)$, $O^2 = (\frac{1}{2},
\frac{2}{3},-p-1)$, $\ O^1 = (-p-\frac{1}{2})$, and
$O^k=\left(-p+\frac{1}{2k}\right)$ or
$\left(-p+\frac{1}{2k+2}\right)$ for $k = 0$ or $3$.
\end{cor}

\nibf{Proof:} Note that the right-handed products, when $P=(0)$, are given when
$\otf=(\frac{1}{2}, \frac{2}{3},-1)$ and
$\onotf=\left(-\frac{1}{2k+2}\right)$ or $\onotf=\left(-\frac{1}{2k}\right)$.
Therefore, the handedness of the product $b(+5,1)$, as
determined through electron microscopy by Crisona \textit{et al.}
\cite{Cri} means that $p=(-p)$, $R = (-1-p)$ and
$O^2 = (\frac{1}{2},\frac{2}{3},p-1) = (\frac{1}{2},
\frac{-1}{3})$, for $p\in\{0,1\}$. By including the
negative sign in $p$ we get the resulting $P$, $R$ and $\otf$.

In the case $k=1$, further electron micrographs of Crisona
\textit{et al.} show the majority of the product is $b(+3,1)$.
Since $R=(-p-1)$, the corresponding right-handed solution is
$O^1=\left(-p-\frac{1}{2}\right)$.

No chirality information
exists for $k=0$ or $3$.
Therefore there are still two possibilities
for $k=0$ or $3$.\hfill $\Box$


\section{Solutions}
\label{Solns}

In this section, we present all tangle solutions for an
integrase acting on DNA with direct or inverted sites.  In
particular, we determine what $R$ must be in each case.

We should note that when $O^k=\oc+\of$ is known to be rational
then the following theorem allows us to find a values of $O^k$ and
$R$ that will satisfy tangle equations arising from an unknotted
substrate and four-plat products.

\begin{thm}[\cite{Dar} Lemma 14, \cite{ES1}]\label{ES-thm} With $c'd-cd'=1:$
\[
N\left(\frac{a}{b}+\frac{c}{d}\right)=N\left(\frac{ad+bc}{ad'+bc'}\right)
\]

\end{thm}

For instance, in our case we have the unknot as substrate, and
torus knots and links (that are the numerator closures of integral
tangles) as products.

\begin{thm}\label{R1}
Given $O^n$, $P = \frac{p}{q}$ and $R$ are rational tangles such
that:
\begin{tabbing}
the\= \kill
\>$N\left(O^n+P\right)=N\left(\frac{1}{\mp n}\right)$, the unknot\\
\>$N\left(O^n+R\right)=N\left(\mp \frac{n-t}{1}\right)$, the
 torus knot or link $T_{\mp(n-t),2}$
\end{tabbing}
Then $O^n=\frac{r\pm pn}{s\mp qn}$ \& $R=\frac{r\pm pt}{-s\pm qt}$
for $rq+ps=1$, for any constant $t$.
\label{findingsols}
\end{thm}

\nibf{Proof.} We begin by noting:

\niit{Claim.} Suppose $p$ and $q$ are relatively prime, and
$rq+ps=1$. Then $r'q+ps'=1$ if and only if $r'=r\pm pt$ and
$s'=s\mp qt$.

\niit{Proof of Claim.} If $r'=r\pm pt$ and $s'=s\mp qt$ then
checking that $r'q+ps'=1$ is a simple calculation. Conversely,
$r'q+ps'=rq+ps$, so $(r'-r)q=p(s-s')$. Since $p$ and $q$ are
relatively prime, then $q|(s-s')$, so $s-s'=qt$ for some $t$.
Therefore $r'q=rq+p(s-s')=rq+pqt$ so $r'=r+pt$. Therefore
$rq+pqt+ps'=1$, \textit{i.e.}, $rq+p(qt+s')=1$, and so $s=s'+qt$,
and hence $s'=s-qt$.

Applying Theorem~\ref{ES-thm} to $N\left(\frac{r\pm pn}{s\mp
qn}+\frac{p}{q}\right)$ and $N\left(\frac{r\pm pn}{s\mp
qn}+\frac{r\pm pt}{-s\pm qt}\right)$ now gives the result: In the
first case, $N\left(\frac{r\pm pn}{s\mp qn}+\frac{p}{q}\right)$\,
$=N\left(\frac{1}{\mp n}\right)$. In the second case, note that
$-p(-s\pm qt)-(-q)(r\pm pt)=1$, so $N\left(\frac{r\pm pn}{s\mp
qn}+\frac{r\pm pt}{-s\pm qt}\right)$
$=N\left(\pm\frac{t-n}{-1}\right)$. \hfill$\Box$

\smallskip

Since for inverted repeats
$N\left(O^k+R\right)=T\left(\pm(2k+1),2\right)$, setting
$n=2k+1+t$ gives a solution for the inverted system. For direct
repeats $N\left(O^k+R\right)=T\left(\pm{2k},2\right)$, so setting
$n=2k+t$ gives a solution for the direct system.

Let us illustrate this by picking an arbitrary fraction for $P$. For
$P=\left(\frac{11}{7}\right)$ we get
$\frac{11}{7}=1+\frac{1}{1+\frac{1}{1+\frac{1}{3}}}$ so
$P=[3,1,1,1]$. Note $7(-3)+11(2)=1$, so let $r=-3$ and $s=2$. Then
we get:
%
\[O^k=\frac{r\pm p(2k+t+1)}{s\mp q(2k+t+1)}=-\frac{22k\pm 11t\pm 11 -3}{14k\pm 7t\pm 7 +2}\]
and
\[R=\frac{-3\pm 11t}{-2\pm 7t}=\frac{3\mp 11t}{2\mp 7t}\]
are solutions to the system for the given $P$, for any constant
$t$.


\begin{cor}
If $P$ is an integral tangle $(n)$ then $R$ is either the infinity
tangle, an integral tangle, or the sum of a vertical and integral
tangle.~\label{Rvert}
\end{cor}

\nibf{Proof.} Theorem~\ref{findingsols} tells us that, since if
$r=1-n$ and $s=1$ then $r+ns=1$, $R$ must be
\[
\left(\frac{1-n\pm nt}{-1\pm
t}\right)=\left(n\right)+\left(\frac{1}{t\mp 1}\right).
\]
If $t=0$ or $\pm 2$ then $R$ is an integral tangle. If $t=\pm 1$
then $R$ is the infinity tangle. Otherwise, it is the sum of a
vertical and integral tangle. \hfill$\Box$

In~\cite{us} we showed there are three classes of solutions for
$P$, $\of$ and $\oc$ to both systems. The example illustrated above
lies in Class $3$, described below.

\nibf{Class 1.} $P = (\infty)$ and $\oc$ and $\of$ are integral.
In this case Theorem~\ref{ES-thm} tells us that $R$ must be
integral: $\frac{p}{q}=\frac{1}{0}$, so for $s=1$ and any $r$, say
$r=0$ we get $R=\left(\frac{\pm t}{-1}\right)$. In addition,
simply knowing that $O^k$ is rational, then Theorem~\ref{ES-thm}
gives us that $O^k$ must be integral: $O^k=\left(\frac{\pm
(2k+1+t)}{1}\right)$ for the inverted case.


The simplest example is when $R=(0)$. In that case
$\of+\oc=O^k=\pm(2k+1)$ for inverted, or $O^k=(\pm 2k)$ for direct
\cite{me}. We can, however, choose $R$ to be any integral tangle,
which corresponds to varying $t$, and $O^k$ changes accordingly.

\nibf{Class 2.} $P$ and $\oc$ are integral, and $\of$ (and hence
$O^k$) is $(\infty)$ for at most value of $k$, integral for at
most 2 values of $k$---for all $k$ in both the directed and
inverted systems---and strictly rational otherwise. When $\of$ is
strictly rational, it must be either vertical or the sum of a
vertical and an integral tangle~\cite{us}, possibly of different
signs. (Note that if the vertical and the horizontal twists have
opposite sign the tangle is said to be
in ``non-canonical" form.) As shown in
Theorem~\ref{Rvert}, $R$ must be the infinity tangle, integral, or
the sum of a vertical and integral tangle. However, $R$ and $\of$
for some $k$ cannot both be $(\infty)$.
%
%

The simplest example in this class is when $P=(0)$, $R=(\infty)$
and $O^k = (\pm 1/(2k+1))$ for inverted, and $O^0=(\infty)$ and
$O^k = (\pm 1/2k)$ ($k>0$) for direct, as first considered in
\cite{Cri}. This is biologically equivalent to the first example,
in that the tangle surgery consists in interchanging the two
tangles with no crossings. In general two solutions are
\emph{biologically equivalent} if their three-dimensional
arrangement is the same, but have different projections. See our
previous paper \cite{us} and \cite{Vaz2} for further discussion of
biological equivalence.

Another example, which is biologically non-equivalent, is when
$P=(\pm 2)$. In this case $R=(\pm 1)$ and $O^0$ is infinity for
inverted repeats and $O^0$ and $O^1$ are integral for direct
repeats; all other $O^k$'s are strictly rational. For direct
repeats $O^0=(\mp 1)$, $O^1=(\mp 3)$ or $(\mp \frac{5}{3})$,
$O^2=(\mp \frac{9}{5})$, and $O^3=(\mp\frac{13}{7})$.
For inverted repeats $O^0=(\infty)$ or $(\mp\frac{3}{2})$,
$O^1=(\mp\frac{7}{4})$, $O^2=(\mp\frac{11}{6})$ and
$O^3=(\mp\frac{15}{8})$.

If a solution has $O^k=(\infty)$, for some $k$, as in the second
example, then in fact it must satisfy several conditions:
%
%

\begin{thm}
\label{ThirdSol} Given a solution from Class 2 then:
\begin{enumerate}
\item[($i$)] If $\exists$ $i$ s.t. $O^i_f=(\infty)$, then
$R$ must be vertical or $(\pm 1)$.
\item[($ii$)] If $\exists$ $i$ and $j$ s.t. $O^i_f = (\infty)$
and $O^j_f$ is integral, then $P\in\{0,\pm 2\}$.
\item[($iii$)] If $\exists$ $i$ s.t. $O^i_f=(\infty)$ and $P \neq (0)$,
then $i=0$. If in addition $O^j_f$ is integral, then $j=1$.

\end{enumerate}
\end{thm}
%

\nibf{Proof.} The first two items were shown in~\cite{us}.
The third item can be proved by looking at $O^k$
and $R$, which can be re-written as $R=\left(\frac{1-p(s\mp
t)}{-s\pm t}\right)$ and $O^k=\left(-p+\frac{1}{s\mp
t\mp(2k+1)}\right)$ for inverted, and $O^k=\left(-p+\frac{1}{s\mp
t\mp(2k)}\right)$ for direct, where $P=(p)\neq (0)$. As $R$ is
vertical, this means that either $s\mp t = 0$ or $p(s\mp t)=2$. If
$s \mp t = 0$, then, if $O^i=(\infty)$, that $2i+1=0$ for
inverted, or $2i=0$ for direct. The first is impossible, so $i=0$
for direct. If further, $O^j$ is integral, then $\pm(2j+1)=\pm 1$
for inverted, or $\pm 2j=\pm 1$ for direct. The first case gives
$j=0$, the second is impossible. The second case is similar: Since
in this case $p=\pm 2$, this means $s\mp t=\pm 1$. So $2i = \mp 1$
for direct, which is impossible, and $2i+1= \mp 1$ for inverted,
so $i=0$. If in addition $O^j$ is integral, then
$\pm1\mp(2j+1)=\pm 1$ for inverted, which is impossible, and $\pm
1\mp(2j)=\pm 1$, so $j=1$, or $j=0$ which cannot occur since
$i=0$, so $j=1$.\hfill $\Box$

\smallskip

\nibf{Class 3.} $P$ and $\of$ are strictly rational and and $\oc$
is integral. In this case $R$ is also strictly rational: given
$O^k=\left(\frac{a}{b}\right)$, with $b\neq\pm 1$, if $R=(c)$ then
$N(O^k+R)=N\left(\frac{a+bc}{b}\right)
=N\left(\frac{2k+1}{1}\right)$, which implies $b=\pm 1$, a
contradiction.
%
%

%
%

An illustration of the last class of solutions was given above,
following Theorem~\ref{findingsols}. This class includes many
different (biologically non-equivalent) possible actions.

\smallskip

\nibf{Class 4 (inverted case only).}
%
%
In the preceding section we demonstrated that there is an
additional solution for inverted sites: $R, \ \oc$ and $P = (p)$
are integral, and $\otf$ is Montesinos with two non-integral
tangles, and in fact is the minimal prime tangle \cite{Ble}. In
particular, $P = \pm(p)$, $R = \pm(1+p)$ and $\oc$ are integral,
$\otf$ is the prime Montesinos tangle
$\mp(\frac{1}{2},\frac{2}{3},p-1)$ and $\onotf =
\pm(-p+\frac{1}{2k})$ or $\mp(p+\frac{1}{2k+2})$, where
$p\in\{0,1\}$. This is the only case in which one of the
constituent tangles is not rational, and in this case it is the
smallest prime tangle (plus possibly a single horizontal twist).

\begin{figure}[htb]
\begin{center}
\psfig{file=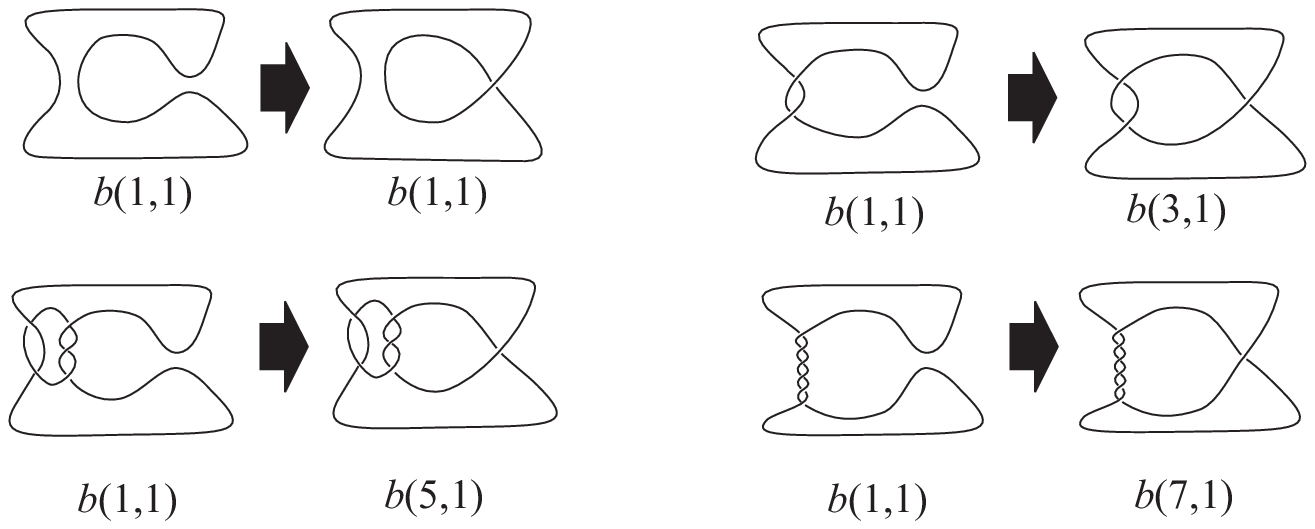,width=3.5in} \caption{The fourth
case for inverted repeats: \mbox{$P=(0), \ R=(-1), \
O^0=(\infty),\ O^1=(-\frac{1}{2}),\ O^2=(\frac{1}{2},
-\frac{1}{3}),\ O^3=(-\frac{1}{6})$}}\label{f:fourthsol}
\end{center}
\end{figure}


\section{Conclusion}

We have completely classified all possible tangle solutions
corresponding to a generic member of the integrase family of
recombinases acting on DNA with either direct or inverted sites.
We proved there are exactly three classes of solutions in the
direct system, and exactly four in the inverted.

The biological assumption is that the protein mechanism is
consistent for a given protein across a variety of substrates. The
integrase family of recombinases is forgiving of
substrates---these members will perform recombination between two
sites on the same (circular, as discussed here, or linear) or
different (circular and/or linear) molecules. However, if an
integrase requires a particular alignment of the two sites when
they are brought together, then a substrate with inverted sites
may have a different number of crossings introduced into $O$ than
a substrate with direct sites. (For example, it was shown that Flp
aligns the sites in antiparallel, and so in the simplest case,
$\oc = (1)$ for inverted and $(0)$ for direct \cite{me}.)
Biologically, the unequal number of solutions between the inverted
and direct cases means that the fourth case is not likely, as this
would mean that $O^k$ differs not only in the number of crossings
but in the actual type of tangle (prime versus rational). This
would in turn imply a different mechanism for bringing the sites
together depending whether they are in direct or inverted repeats
(only for one product, corresponding to $k=2$). A small protein
being able to ``sense'' the global orientation (inverted versus
direct) of the two sites would be unusual. However, there are
examples of proteins which act in response to global
characteristics, notably Topoisomerase II, which preferentially
changes crossings to \textit{un}knot DNA~\cite{Ryben}.

%
%


\subsection{Future Directions}
We conclude with several observations, and avenues for future
biological and mathematical research.

Electron microscopy images by Crisona and colleagues determine
that for Flp in the cases $k=1$ and $k=2$, the products are in
fact torus knots~\cite{Cri}. For the related recombinase,
$\lambda$ integrase acting on LR or PB sites, electron microscopy
has shown that the products are almost exclusively positive torus
knots (inverted) and negative torus links
(direct)~\cite{Cri},~\cite{Speng}. Biologically, one would expect
a fixed, precise mechanism for this family of recombinases, which
would predict that for Flp and other integrases, every product
should be a torus knot or link. However, although one can
experimentally demonstrate that higher-crossing products exist and
determine their respective crossing number, their precise
knot/link type has not been experimentally verified for any other
integrase family member. As the crossing number increases,
resolving the precise knot type becomes more crucial (and
difficult) since the number of knots with a given crossing number
increases dramatically---for instance, there are 1,701,936 knots
with $\leq$ 16 crossings \cite{HosThis}.

\nibf{Question 1:} \emph{For Flp and other integrases (not
$\lambda$ Int), can we obtain experimental confirmation that every
product is of the form $b(2k,1)$ (direct) or $b(2k+1,1)$
(inverted)}?

If they are, one might be able develop similar Dehn surgery
arguments such as the result from~\cite{KMOS} used here to
restrict $O^k$. For example, for $k \leq 10$, we can eliminate the
possibility of $O^k$ having as double branch cover a satellite
knot complement by work of Bleiler and Litherland~\cite{BL}, and,
in many cases, a torus knot complement by Moser and Gordon
\cite{Mos,Gor}. So in many cases, the first step would be to
consider $O^k$ such that $\wt{O^k}$ is a hyperbolic knot
complement.

\smallskip

Second, we consider the chirality of the resulting products.
Crisona \textit{et al.} have characterized the first two nontrivial
products of Flp-mediated inversion as almost exclusively positive
torus knots \cite{Cri}.

\nibf{Question 2}: \emph{Can we determine the chirality of any of
the Flp deletion products, or the higher crossing inversion
products?}

If so, the classes of solutions could be tightened signifcantly by
removing the possibility of mirror images, as in Corollary
\ref{chiral}.

\smallskip

Third, although the third class of solutions is mathematically
possible, biological considerations such as DNA's stiffness
impeding a high number of crossings, make them biologically
unlikely. For this reason, biologists often assume in the tangle
model that $P$ and $R$ each have 0 or 1 crossings.

Other work, notably~\cite{Cri} and~\cite{SECS}, have incorporated
a number of biologically reasonable assumptions into the tangle
model which reduce the number of putative tangle solutions. For
example the Generalized Random Collision assumes that the original
DNA is exclusively negatively supercoiled (\textit{i.e.}, $\of =
(n)$ or $(-\frac{1}{n})$ for $n\geq 0$), which then in turn
restricts $\oc$ and $P$~\cite{Cri}, \cite{SECS}. These assumptions
are biologically reasonable, and preclude tangle solutions such as
those of Class 3. But topologically, the more exotic cases cannot
be excluded.

%

\nibf{Question 3:} \emph{Can the third class of solutions be
eliminated or restricted experimentally, or computationally using
models (\textit{e.g.},~\cite{MB}) of DNA flexibility?}

If the number of crossings of the constituent tangles can be
bounded with any degree of certainty, this reduces the solutions
from an infinite number to a mere handful. For example, we
consider the serine recombinase Xer acting on circular DNA with
direct \textit{psi} sites, whose corresponding tangle equations
are $N(O+P)= b(1,1)$ and $N(O+R)=N\left(\frac{4}{1}\right) =
b(4,1)$~\cite{Colloms}. Assuming all tangles are rational and
$P=(0)$, and so $O=(\frac{1}{r})$, then Darcy showed that
$R\in\{\frac{1}{j},\frac{3}{3+j},\frac{5}{5+j},\frac{4k-1}{4+j(4k-1)}\}$~\cite{Dar}.
Additionally, Vazquez \emph{et al}, by carefully analyzing the
biological data, and assuming any nontrivial topology is in $\oc$,
show $P$ is trivial (and hence in their setting can be chosen to
be $(0)$), and show $O$ is rational using techniques similar to
\cite{us}. Further they assume $R$ is integral or $\infty$, and
show there exist 3 solutions: $O = (\frac{-1}{3}), R = (-1)$ or $O
= (\frac{-1}{5}), R = (+1)$ or $O = (\frac{-1}{4}), R = (\infty)$,
which can be biologically equivalent \cite{Vaz2}.

Additionally given torus link substrates constructed by $\lambda$
Int, Xer-mediated recombination results in more complex knots:
$N(O_2+P)=N\left(\frac{6}{1}\right)$, and $N(O_2+R)=$
seven-crossing product (knot), and
$N(O_3+P)=N\left(\frac{8}{1}\right)$ $N(O_3+R)=$ nine-crossing
product~\cite{Bath}. Darcy also uses this information to analyze
how differing seven-crossing products reduce the number of
putative solutions~\cite{Dar}.

\nibf{Question 4} \emph{Can the precise type of knot resulting in
the action of Xer on $N\left(\frac{6}{1}\right)$ be determined?
Can similar techniques to those used here be used to examine the
Xer system in their full generality, without simplifying
assumptions?}

\smallskip

Finally, we reiterate that, although we have motivated our work by
considering Flp, our results can be applied to any protein in the
integrase family, such as Cre or $\lambda$ Int (acting on LR
sites), whose products are the above torus knots or links.

\medskip

\nibf{Acknowledgements.} We wish to thank Francis Bonahon, Cameron
Gordon and John Luecke for a number of illuminating discussions.
We also would like to thank Makkuni Jayaram for introducing us to
Flp. Our appreciation also goes to Isabel Darcy for her careful
reading of an earlier draft, particularly for pointing out that
Montesinos tangles must be considered in their generalized form,
and suggesting Corollary \ref{newcor}.

DB would also like to thank the Institut des Hautes \'{E}tudes
Scientifiques for its hospitality during the writing of this work.

DB was supported by a grant from the National Science Foundation's
Division of Mathematical Sciences. CVM was supported in part by
the Presidential Excellence Summer Scholarly Activity Grant from
St.\ Edward's University.


\bibliographystyle{amsplain}



\end{document}